\documentclass[reqno]{amsart}

\usepackage{amsmath}
\usepackage{amsfonts}
\usepackage{amssymb}
\usepackage{verbatim}

\newtheorem{theorem}{Theorem}[section]

\newtheorem{corollary}[theorem]{Corollary}

\theoremstyle{definition}
\newtheorem{definition}[theorem]{Definition}

\theoremstyle{remark}

\numberwithin{equation}{section}

\begin{document}
\title[Examples of noncommutative manifolds]{Examples of noncommutative manifolds:\\ complex tori and spherical manifolds}
\author{Jorge Plazas}
\address{ Max Planck Institute for Mathematics,
Vivatsgasse 7,
Bonn 53111,
Germany}
\email{plazas@mpim-bonn.mpg.de}
\subjclass[2000]{Primary 58B34, 46L87; Secondary 46L89, 16S38}
\keywords{noncommutative tori, noncommutative spherical manifolds}

\begin{abstract}
We survey some aspects of the theory of noncommutative manifolds 
focusing on the noncommutative analogs of two-dimensional tori and 
low-dimensional spheres.
We are particularly interested in those aspects of the theory that link 
the differential geometry and the algebraic geometry of these spaces.
\end{abstract}
\maketitle

\setcounter{section}{-1}
\section{Introduction}
\label{intro}

These notes are based on a series of lectures given at the International Workshop on Noncommutative Geometry held at Tehran in September 2005. The aim of the lectures was to introduce some examples of noncommutative manifolds as well as some of the tools and notions needed for their study. We focused on those examples in which recent developments have provided links between the smooth and the algebraic aspects of noncommutative  geometry.    

In the first part of the lectures some aspects of the theory of noncommutative two-dimensional tori are treated. Noncommutative tori are standard prototypes of noncommutative spaces. Since the early stages of noncommutative geometry these spaces have been central examples arising naturally in various contexts (cf. \cite{Connes1, rieffel1}). Their topology and differential geometry are by now well understood. In particular, vector bundles on these spaces and the corresponding theory of Morita equivalences can be characterized. Connections on these vector bundles arising as liftings of the natural derivations on the noncommutative torus give rise to a rich theory  (cf. \cite{Connes1, ConnesRieffel}). Recently this theory has been recast in the context of complex algebraic geometry. The study of categories of holomorphic bundles has thrown some light on the underlying algebraic structures of these spaces 
(cf. \cite{DiS, Po, PoS, Schwarz1, yo}). Noncommutative tori can be viewed as noncommutative limits of elliptic curves. Moreover, it is believed that noncommutative tori may play a role in number theory analogous to the role played by elliptic curves 
(cf. \cite{Manin1,Manin2}). This opens a new perspective for the fruitful study of these spaces. We begin by discussing the differential topology of noncommutative tori. We then discuss how the theory of vector bundles on noncommutative tori and the corresponding $K$-theory lead to the existence of a class of noncommutative tori whose structure is related to arithmetic data. Finally we study holomorphic connexions on these spaces and the homogeneous coordinate rings arising from them.

The second class of examples we treat are noncommutative spaces defined in terms of a 
$K$-theoretic equation whose solutions in the commutative case correspond to the classical spheres; we call a noncommutative space of this type a spherical manifold. The characterization of the geometry of spheres from a spectral point of view raised the question about the existence of such examples (cf. \cite{Connes4}). The first noncommutative examples of spherical manifolds were constructed as a one-parameter family of deformations of the algebra of smooth functions on the sphere $S^4$ 
(cf. \cite{LandiConnes}). These noncommutative spheres are a part of a broader class of noncommutative manifolds known as $\theta$-deformations. Shortly afterwards  
the case of spherical manifolds in dimension 3 was carried out systematically 
(cf. \cite{CoDVI,CoDVII}). For this purpose new tools have to be introduced; these tools provide links between the algebraic and the differential aspects of noncommutative geometry.

I am happy to thank the organizers of the conference and the Institute for Studies in Theoretical Physics and Mathematics at Tehran for giving me the opportunity to take part in this very nice conference. I am also grateful to Arthur Greenspoon, Masoud Khalkhali and Matilde Marcolli; their comments helped me to greatly improve the present manuscript. 

\section{Noncommutative tori}

In what follows we will denote by $\mathbb{T}^{2}$ the two-dimensional torus 
$\mathbb{T} \times \mathbb{T}$, where $ \mathbb{T} = S^1=\mathbb{R}/\mathbb{Z}$.
We can consider different structures on $\mathbb{T}^{2}$; these different structures will give rise to algebras of functions of various levels of regularity. For instance, considering  $\mathbb{T}^{2}$ as a topological space leads to the study of its algebra of continuous functions and considering  $\mathbb{T}^{2}$ as a smooth manifold leads to the study of its algebra of smooth functions. Once we endow the smooth manifold $\mathbb{T}^{2}$ with a complex structure it becomes a Riemann surface of genus one, i.e. an elliptic curve over $\mathbb{C}$. Once enriched with this algebro-geometric structure $\mathbb{T}^{2}$ plays a fundamental role in number theory.

In this section we will study the noncommutative analogs of $\mathbb{T}^{2}$. These noncommutative tori are defined in terms of noncommutative algebras which can be obtained as deformations of some algebras of functions on $\mathbb{T}^{2}$ and share many of their structural properties. In particular it makes sense to study the geometric behavior of these spaces in a framework which generalizes the classical setting. One advantage of passing to the noncommutative setting comes from the fact that in many cases the presence of noncommutativity enriches the classical picture. For instance, as will be seen below, noncommutative tori exhibit already at the topological and smooth levels an arithmetic behavior close to the arithmetic behavior of elliptic curves.

\subsection{Noncommutative tori as topological spaces}

As a topological space $\mathbb{T}^{2}$ its characterized by $C(\mathbb{T}^2)$, the algebra of continuous functions from  $\mathbb{T}^2$ to the complex numbers. 
In general the algebra of continuous complex-valued functions on any compact Hausdorff space is a unital commutative $C^{*}$-algebra. Conversely, by a theorem of Gelfand, any unital commutative $C^{*}$-algebra can be realized as the algebra of continuous functions on some compact Hausdorff space. This is one of the departure points of noncommutative geometry. We can drop the commutativity condition and consider a noncommutative $C^{*}$-algebra as the algebra of continuous functions on some noncommutative topological space defined by duality with this algebra. In this section we will describe the noncommutative analogs of $C(\mathbb{T}^2)$. The reader may refer to \cite{Davidson}, \cite{GraciaVarilly} or \cite{Masoud} for the basic concepts and definitions of the theory of $C^{*}$-algebras.

Lets start by taking a closer look at $C(\mathbb{T}^{2})$. The algebra $C(\mathbb{T}^{2})$ can be realized as the universal $C^{*}$-algebra generated by two commuting unitaries $U$ and $V$. One may for instance take $U$ and $V$ to be the functions $e^{2 \pi \imath \varphi_1}$ and $e^{2 \pi \imath \varphi_2}$ where 
$\varphi_{j} \in \mathbb{R}/\mathbb{Z}, \; j=1,2$ give the coordinate functions on 
$\mathbb{T}^{2}$. Any element of $C(\mathbb{T}^{2})$ admits a Fourier expansion in terms of powers of these unitaries. Moreover, the algebra $C(\mathbb{T}^{2})$ maps to any other $C^{*}$-algebras generated by two commuting unitaries. Therefore $C(\mathbb{T}^{2})$ is the \emph{universal $C^{*}$-algebra generated by two commuting unitaries $U$ and $V$}.

At a topological level noncommutative tori are defined by their algebras of continuous functions, which are noncommutative deformations of $C(\mathbb{T}^2)$ defined in terms of a similar universal property where we have dropped the commutativity condition on the generating unitaries:

\begin{definition}
Given  $\theta \in \mathbb{R}$ we define $A_{\theta} = C(\mathbb{T}^2_\theta)$, the {\it algebra of continuous  functions on the noncommutative torus $\mathbb{T}^2_\theta$}, as the universal $C^{*}$-algebra generated by two unitaries $U$ and $V$ subject to the relation:
\begin{eqnarray}
\label{commrel}
UV=e^{2\pi \imath \theta}VU.
\end{eqnarray}
\end{definition}

\noindent Note that if we take $\theta = 0 $ we recover $C(\mathbb{T}^2)$.

Fix $\theta \in \mathbb{R}$. We can use a standard procedure to construct the $C^{*}$-algebra 
$A_{\theta}$ starting from the unitaries $U$ and $V$ together with the relation (\ref{commrel}). 
First we construct a $*$-algebra $C_{alg}(\mathbb{T}^2_\theta)$ by taking all $\mathbb{C}$-linear combinations of products of the generators $U,\, U^*, \, V$ and $V^*$ and imposing the relations:
\begin{eqnarray*}
UU^* &=& U^*U \,=1 \\
VV^* &=& V^*V  \,= 1 \\
UV   &=& e^{2\pi \imath \theta}VU.
\end{eqnarray*}
Equivalently $C_{alg}(\mathbb{T}^2_\theta)$ is the quotient of the free algebra 
$\mathbb{C}\langle U,U^*,V,V^* \rangle$ by the above relations. The $*$-algebra structure is the obvious one defined by $*: U\mapsto U^{*}$ and 
$*: V\mapsto V^{*}$. Then we consider representations of  $C_{alg}(\mathbb{T}^2_\theta)$ as an algebra of bounded operators. 
These are given by $*$-morphisms
\begin{eqnarray*}
\rho: C_{alg}(\mathbb{T}^2_\theta) \rightarrow \mathcal{B}(\mathcal{H}_{\rho})
\end{eqnarray*}
where $\mathcal{B}(\mathcal{H}_{\rho})$ is the algebra of bounded operators on a Hilbert space 
$\mathcal{H}_{\rho}$. In this way we may endow $C_{alg}(\mathbb{T}^2_\theta)$ with the norm 
\begin{eqnarray*}
\lVert a \rVert = \sup_{\rho} \lVert \rho(a) \rVert, 
&& a \in C_{alg}(\mathbb{T}^2_\theta)
\end{eqnarray*}
where the supremum is taken over all the representations of 
$C_{alg}(\mathbb{T}^2_\theta)$ and the norm on the right is the operator norm. This supremum will be finite since unitary operators have norm one and the images of the generators $U$, $U^*$, $V$ and $V^*$ are all unitary. The $C^{*}$-completion of 
$C_{alg}(\mathbb{T}^2_\theta)$ under this norm is the algebra 
$A_\theta = C(\mathbb{T}^2_\theta)$.

For the above construction to work we have to show the existence of nontrivial representations of $C_{alg}(\mathbb{T}^2_\theta)$. For this purpose we consider the unitary operators $u$ and $v$ on $L^2(S^1)$ given by
\begin{eqnarray}
\label{rep1}
(uf)(t) &=& e^{2 \pi \imath \theta} f(t) \\
(vf)(t) &=& f(x- \theta) 
\end{eqnarray}
Since $vu=e^{2\pi \imath \theta}uv$ the map $U\mapsto u$, $V\mapsto v$ gives a representation of $C_{alg}(\mathbb{T}^2_\theta)$ as a subalgebra of the algebra of bounded operators on $L^2(S^1)$.

From the definition we see that if $\theta=\theta' + n$ for some integer $n$ then the corresponding algebras $A_\theta$ and $A_{\theta'}$ are isomorphic as $C^{*}$-algebras. So whenever it is necessary we can assume $0 \leq \theta < 1$. Moreover, by universality of the construction we obtain an isomorphism between $A_\theta$ and $A_{1-\theta}$ by interchanging the roles of the generators.  Thus we may equally well assume 
that $0 \leq \theta \leq \frac{1}{2}$. It can be shown that for different values of 
$\theta$ in this interval one obtains nonisomorphic $C^{*}$-algebras.

Many of the tools used for the study of topological spaces can be extended to the study of 
$C^{*}$-algebras. One important phenomenon that arises in this context is that topological invariants like $K$-theory and cyclic cohomology are invariant under strong Morita equivalence of $C^{*}$-algebras. This notion coincides in the case of unital $C^{*}$-algebras with the notion of Morita equivalence for rings which identifies rings whose categories of right modules are equivalent (cf. \cite{GraciaVarilly}). In the next sections we will look more closely at Morita equivalences between noncommutative tori. In the context of $C^{*}$-algebras Morita equivalences between noncommutative tori are characterized by the following result:  

\begin{theorem}\emph{(\cite{rieffel1})}
\label{Morrita}
Let $SL_2 (\mathbb{Z})$ act on $\mathbb{R}$ by fractional linear transformations; thus
\begin{eqnarray*}
g \theta &=& \frac{a\theta +b}{c\theta + d} 
\end{eqnarray*}
where
\begin{eqnarray*}
g= \left(
\begin{array}{cc}
  a & b \\
  c & d \\
\end{array}
\right)\in SL_2 (\mathbb{Z}), \qquad \theta \in \mathbb{R}. 
\end{eqnarray*}
Then the $C^{*}$-algebras $A_\theta$ and $A_{\theta'}$ are strongly Morita equivalent if and only if there exist a matrix $g \in SL_2 (\mathbb{Z}) $ such that $\theta' =g \theta$
\end{theorem}

Taking into account that any two rational numbers are conjugate to each other by a fractional linear transformation in $SL_2(\mathbb{Z})$ we obtain the following corollary:  
\begin{corollary}
If $\theta \in \mathbb{Q}$ then the $C^{*}$-algebra $A_\theta$ is strongly Morita equivalent to  $A_{0}= C(\mathbb{T}^2)$. 
\end{corollary}

In what follows we will assume unless otherwise stated that $\theta$ is an irrational number.

The compact group $\mathbb{T}^2$ acts on the algebra $A_\theta$. This action can be given in terms of the generators $U$ and $V$ by:
\begin{eqnarray*}
\alpha_{\varphi}(U) &=& e^{2 \pi \imath \varphi_1}U  \\
\alpha_{\varphi}(V) &=& e^{2 \pi \imath \varphi_2}V
\end{eqnarray*}
where $ \varphi = (\varphi_1,\varphi_2) \in \mathbb{T}^2$.
This action is continuous in the sense that for  any $a\in A_\theta$ the function 
\begin{eqnarray}
\nonumber
\mathbb{T}^2 & \rightarrow & A_\theta \\
\label{act}
\varphi &\mapsto& \alpha_{\varphi}(a)
\end{eqnarray} is continuous. 

By continuity of the above action the following integral makes sense as a limit of Riemann sums for any $a\in C(\mathbb{T}^2_\theta)$: 
\begin{eqnarray*}
\label{laint}
\int_{\mathbb{T}^2} \alpha_{\varphi} (a) \, d\varphi 
\end{eqnarray*}
Denote by $\mathbf{1}_{A_\theta}$ the unit of the algebra $A_\theta$. It can be shown that this integral takes values in $\mathbb{C} \mathbf{1}_{A_\theta}$ and thus it defines a linear functional which once normalized induces a trace 
\begin{eqnarray}
\label{trace1}
\chi: A_\theta \rightarrow \mathbb{C}
\end{eqnarray}

\begin{theorem}\emph{(\cite{Effros})}
Let $\theta$ be an irrational number. Then 
$\chi$ is the unique normalized trace on $A_{\theta}$ invariant under the action of $\mathbb{T}^2$.
\end{theorem}

This fact has the following important consequence (cf. \cite{Davidson})
\begin{theorem}\emph{(\cite{Effros})}
Let $\theta$ be an irrational number. Then 
$A_{\theta}$ is simple, i.e. it has no nontrivial two-sided ideals.
\end{theorem}

The $C^{*}$-algebras $A_\theta$ arise naturally in various contexts. Let us finish this section by briefly mentioning two of them.

It is possible to associate a crossed product  $C^{*}$-algebra to an action of a discrete group on a topological space. Let $\theta$ be irrational and consider the action $\tilde{\theta}$ of $\mathbb{Z}$ on $S^{1}$ generated by the rotation by an angle $\theta$. The corresponding crossed product algebra is then isomorphic to $A_\theta$. Because of this fact the algebras $A_\theta$ are referred to as \emph{irrational rotation algebras} 
(cf. \cite{rieffel1}).

It is also possible to associate a $C^{*}$-algebra to a smooth foliation on a manifold. Consider the Kronecker foliation $d \varphi_1= \theta d \varphi_2$ on $\mathbb{T}^{2}$. The $C^{*}$-algebra associated to this foliation is isomorphic to $A_\theta$ 
(cf. \cite{Connes1}). 

\subsection{Smooth functions on noncommutative tori}

In this section we study the noncommutative analogs of the algebra $C^\infty (\mathbb{T}^{2})$ consisting of smooth complex-valued functions on $\mathbb{T}^{2}$. As before we start by taking a closer look at the commutative case. In this case any element of $C(\mathbb{T}^{2})$ admits a Fourier expansion and smooth functions are characterized among continuous functions as those functions whose coefficients in the corresponding Fourier expansion decay rapidly at infinity. In view of this fact we make the following definition:

\begin{definition}
Given  $\theta \in \mathbb{R}$ we define $\mathcal{A}_\theta$, the {\it algebra of smooth functions on the noncommutative torus $\mathbb{T}_\theta^{2}$}, as the algebra of formal power series in two unitaries $U$ and $V$ with rapidly decreasing coefficients and multiplication given by the relation $UV=e^{2\pi \imath \theta}VU$:
\begin{eqnarray*}
\mathcal{A}_\theta &=& C^\infty (\mathbb{T}_\theta^{2}) \\
 &=& \{ a =
\sum_{n,m \in \mathbb{Z}} a_{n,m}U^nV^m  \, | \, \{a_{n,m} \}
\in \mathcal{S}(\mathbb{Z}^2) \} \nonumber
\end{eqnarray*}
\end{definition}

The algebra $\mathcal{A}_\theta$ is a pre-$C^{*}$-algebra whose $C^{*}$-completion is isomorphic to $A_\theta$. As a subalgebra $\mathcal{A}_\theta$  can be obtained as the set of 
\emph{smooth elements of $A_\theta$ for the action of $\mathbb{T}^{2}$}; that is, elements $a \in A_\theta$ for which the map (\ref{act}) is smooth. 

In the algebra $\mathcal{A}_\theta$ the unique trace $\chi$ defined by (\ref{laint}) is given  by
\begin{eqnarray*}
\chi (\sum a_{n,m}U^nV^m ) = a_{0,0}
\end{eqnarray*}

The action of $\mathbb{T}^2$ on $A_\theta$ induces an action of its Lie algebra $L= \mathbb{R}^2$ on $\mathcal{A}_\theta$ given by the derivations 
\begin{eqnarray} 
\label{derivation1} 
\delta_{1}(U)= 2\pi \imath U; && \delta_{1}(V)= 0 \\
\label{derivation2} 
\delta_{2}(U)= 0; && \delta_{2}(V)= 2\pi \imath V
\end{eqnarray}

These derivations determine the Fr\'{e}chet structure of $\mathcal{A}_\theta$. 
The relation between $\mathcal{A}_\theta$ and $A_\theta$ parallels the relation between
$C^\infty (\mathbb{T}^2)$ and $C(\mathbb{T}^2)$. We refer the reader to the seminal paper 
\cite{Connes1} and the survey \cite{rieffel2} for the main results about the algebras
$A_\theta =C(\mathbb{T}_\theta^2)$ and  $\mathcal{A}_\theta = C^\infty (\mathbb{T}_\theta^2) $.
In the next section we exploit the corresponding relation between smooth vector bundles and continuous vector bundles over the noncommutative torus. 

\subsection{Morita equivalences and real multiplication}

By the Serre-Swan theorem the theory of vector bundles over $\mathbb{T}^2$ is equivalent to
the theory of finite type projective modules over the algebra $C(\mathbb{T}^2)$. To each complex vector bundle over $\mathbb{T}^2$ one associates the $C(\mathbb{T}^2)$-module of its global sections. Smooth bundles correspond to finite type projective modules over $C^\infty (\mathbb{T}^2)$ and every continous vector bundle over $\mathbb{T}^2$ is equivalent to a smooth vector bundle.

We consider projective finite type right $A_\theta$-modules as continuous vector bundles over 
$\mathbb{T}_\theta^2$ and projective finite type right $\mathcal{A}_\theta$-modules as smooth vector bundles over $\mathbb{T}_\theta^2$. If $\tilde{E}$ is a projective finite type right $A_{\theta}$-module then there exists a projective finite type right $\mathcal{A}_\theta$-module $E$ such that one has an isomorphism of right $A_{\theta}$-modules:
\begin{eqnarray*}
\tilde{E} \simeq E \otimes_{\mathcal{A}_\theta} A_{\theta}
\end{eqnarray*}
Therefore, as in the commutative case, the categories of smooth and continuous vector bundles over $\mathbb{T}_\theta^2$ are equivalent (cf. \cite{Connes1}). In what follows we will restrict to $\mathcal{A}_\theta$-modules. 

For any positive integer $n$ the trace $\chi$ can be extended to a trace $\mathrm{Tr}_{\chi}$ on the matrix algebra
\begin{eqnarray*}
\mathcal{M}_n(\mathcal{A}_\theta) &=& \mathcal{A}_\theta \otimes \mathcal{M}_n(\mathbb{C}) \\
&\simeq & \mathrm{End}(\mathcal{A}_\theta^n)
\end{eqnarray*} 
We can use this trace to define a rank function for vector bundles over $\mathbb{T}_\theta^2$. A
right $\mathcal{A}_\theta$-module $E$ is  projective of finite type  if and only if there exists some $n$ and an idempotent $e=e^{2}=e^{*}$ in $\textsl{M}_n(\mathcal{A}_\theta)$ such that $E \simeq e\mathcal{A}_\theta^n$; thus we can define the rank of $E$ by
\begin{eqnarray}
\mathrm{rk}(E)= \mathrm{Tr}_{\chi}(e)
\end{eqnarray}

In what follows unless otherwise stated a right (resp. left) $\mathcal{A}_\theta$-module will always mean a right (resp. left) projective finite type $\mathcal{A}_\theta$-module.

Let $\theta\in\mathbb{R}$ be irrational. Following \cite{Connes1} we define, for any pair $c,d\in \mathbb{Z}$, $c>0$, a right $\mathcal{A}_\theta$-module $E_{d,c}(\theta)$ given by the following action of $\mathcal{A}_\theta $ on the Schwartz space $\mathcal{S} (\mathbb{R} \times \mathbb{Z}/ c\mathbb{Z}) = \mathcal{S}(\mathbb{R})^c$: 
\begin{eqnarray}
(fU)(x,\alpha)&=& f(x- \frac{c\theta+d}{c},\alpha -1) \\
(fV)(x,\alpha)&=& \exp (2 \pi \imath (x- \frac{\alpha d}{c})) f(x,\alpha)
\end{eqnarray}

\begin{theorem}\emph{(\cite{Connes1})}
The rank of $E_{d,c}(\theta)$ is $| c \theta +d |$. If $E$ is any right $\mathcal{A}_\theta$-module with $\mathrm{rk}(E) = | c \theta +d |$ then $E \simeq E_{d,c}(\theta)$. 
\end{theorem}

The $K_{0}$ group of $\mathcal{A}_\theta$, $K_{0}(\mathcal{A}_\theta)$, is by definition the enveloping group of the abelian semigroup given by isomorphism classes of right $\mathcal{A}_\theta$-modules together with direct sum. The rank function $\mathrm{rk}$ extends to an injective morphism 
\begin{eqnarray}
\mathrm{rk}: K_{0}(\mathcal{A}_\theta) \rightarrow \mathbb{R}
\end{eqnarray}
whose image is $\mathbb{Z} \oplus \theta \mathbb{Z}$. Therefore, one gets an ordered structure on $K_{0}(\mathcal{A}_\theta)$ given by the isomorphism 
\begin{eqnarray}
K_{0}(\mathcal{A}_\theta) \simeq \mathbb{Z} \oplus \theta \mathbb{Z} \subset \mathbb{R}
\end{eqnarray}

By a cancellation theorem due to Rieffel the abelian semigroup of isomorphism classes of right $\mathcal{A}_\theta$-modules is a cancellation semigroup (cf. \cite{riff3}).  
This fact together with the fact that $\mathrm{rk}$ is injective imply that right $\mathcal{A}_\theta$-modules are classified up to isomorphism by their rank and that any finite type projective right $\mathcal{A}_\theta$-module is either free or isomorphic to a right module of the form $E_{d,c}(\theta)$.

If $c$ and $d$ are relatively prime we say that $E_{d,c}(\theta)$ is a {\it basic} $\mathcal{A}_\theta $-module. This being the case the pair $d,c$ can be completed to a matrix  
\begin{eqnarray}
\label{lag}
g= \left(
\begin{array}{cc}
  a & b \\
  c & d \\
\end{array}
\right)\in SL_2 (\mathbb{Z})
\end{eqnarray}

We write $E_g(\theta)$  for the module $E_{d,c}(\theta)$. By definition the degree of $E_g(\theta)$ is taken to be $c$. We also define the degree of a matrix $g \in SL_2 (\mathbb{Z})$, given as above, by 
$\mathrm{deg}(g)=c$.

Let $SL_2 (\mathbb{Z})$ act on $\mathbb{R}$ by fractional linear transformations. Let $g\in SL_2 (\mathbb{Z})$ be as above and denote by $U'$ and $V'$ two generating unitaries of the algebra $ \mathcal{A}_{g \theta}$. We can define a left action of the algebra $\mathcal{A}_{g \theta}$ on $E_g$ by: 
\begin{eqnarray}
(U'f)(x,\alpha)&=& f \left(x- \frac{1}{c},\alpha -a \right) \\
(V'f)(x,\alpha)&=& \exp (2 \pi \imath (\frac{x}{c\theta +d}-\frac{\alpha}{c})) f(x,\alpha)
\end{eqnarray}
This action gives an identification: 
\begin{eqnarray}
\mathrm{End}_{\mathcal{A}_\theta}(E_g(\theta)) \simeq \mathcal{A}_{g\theta}.
\end{eqnarray}

The tensor product of basic modules is again a basic module. More precisely, given $g_{1}, g_{2} \in SL_2 (\mathbb{Z})$ there is a well defined pairing of right $\mathcal{A}_\theta$-modules:
\begin{eqnarray}
\label{producto1}
t_{g_{1}, g_{2}}: E_{g_{1}}(g_{2}\theta)\otimes_{\mathbb{C}}  E_{g_{2}} (\theta) \rightarrow 
E_{g_{1} g_{2}}(\theta)
\end{eqnarray}
This map gives rise to an isomorphism of $\mathcal{A}_{g_{1}g_{2}\theta} - \mathcal{A}_{\theta}$ bimodules: 
\begin{eqnarray}
E_{g_{1}}(g_{2}\theta)\otimes_{\mathcal{A}_{g_{2}\theta}}  E_{g_{2}} (\theta) \rightarrow 
E_{g_{1} g_{2}}(\theta).
\end{eqnarray}
In particular, if $g\theta = \theta$ one has an isomorphism 
\begin{eqnarray}
\underbrace{E_{g}(\theta)\otimes_{\mathcal{A}_ {\theta}} \dots
\otimes_{\mathcal{A}_ {\theta}}  E_{g} (\theta)}_{n} 
\simeq  E_{g^{n}} (\theta).
\end{eqnarray}

Taking into account the remarks at the beginning of this section  we can now go back to the setting of  Section~\ref{Morrita}. We say that two noncommutative tori $\mathbb{T}_{\theta'}^2$ and $\mathbb{T}_{\theta}^2$ are Morita equivalent if their corresponding algebras of smooth functions are Morita equivalent. Thus 
$\mathbb{T}_{\theta'}^2$ and $\mathbb{T}_{\theta}^2$ are Morita equivalent if and only if there exist an $\mathcal{A}_{\theta'}$-$\mathcal{A}_{\theta}$-bimodule which is projective and of finite type both as a left $\mathcal{A}_{\theta'}$-module and as a right $\mathcal{A}_{\theta}$-module. We will consider the category whose objects are noncommutative tori and whose morphisms are given by isomorphism classes of finite type projective bimodules over the corresponding algebras of smooth functions. Composition is provided by tensor product over the corresponding algebra. An isomorphism in this category is called a \emph{Morita equivalence}. From the discussion above we see that given a real number $\theta$ and a matrix $g \in SL_2 (\mathbb{Z})$ the noncommutative torus
$\mathbb{T}_{g \theta}^2$ and $\mathbb{T}_{\theta}^2$ are Morita equivalent. The inverse of the morphism represented by the $\mathcal{A}_{g \theta}$-$\mathcal{A}_{\theta}$-bimodule $E_{g} (\theta)$ is the morphism represented by the $\mathcal{A}_{\theta}$-$\mathcal{A}_{g \theta}$-bimodule$E_{g^{-1}} (g \theta)$.  By a result of Rieffel these are the only possible Morita equivalences in the category of noncommutative tori (cf. \cite{rieffel1}). More precisely, two noncommutative tori $\mathbb{T}_{\theta'}^2$ and $\mathbb{T}_{\theta}$ are Morita equivalent if and only if there exists $g \in SL_2 (\mathbb{Z})$ such that $\theta' =g \theta$.

If $g\theta =\theta$ then $E_g(\theta)$ has the structure of $\mathcal{A}_\theta$-bimodule. An irrational number $\theta\in\mathbb{R}\setminus \mathbb{Q}$ is a fixed point of a fractional  linear transformation $g \in SL_2 (\mathbb{Z})$ if and only if it generates a quadratic extension of $\mathbb{Q}$. 

\begin{definition}
The noncommutative torus $\mathbb{T}_\theta^2$ with algebra of smooth functions 
$\mathcal{A}_\theta$ is a {\it real multiplication noncommutative torus} if the parameter $\theta$ generates a quadratic extension of $\mathbb{Q}$.
\end{definition}

Thus $\mathbb{T}_\theta^2$ is a  real multiplication noncommutative torus if and only if it has nontrivial Morita autoequivalences.

In \cite{Manin1} Manin proposed the use of noncommutative tori as a geometric framework under which to attack the explicit class field theory problem for real quadratic extensions of $\mathbb{Q}$. The explicit class field theory problem for number fields asks for explicit generators of the maximal abelian extension of a given number field and the corresponding Galois action of the abelianization of the absolute Galois group on these generators. The only number fields for which a complete solution of this problem is known are the imaginary quadratic extensions of $\mathbb{Q}$ and $\mathbb{Q}$ itself. Elliptic curves whose endomorphism ring is nontrivial correspond to lattices generated by imaginary quadratic irrationalities and play a central role in the solution of the explicit class field theory problem for the corresponding imaginary quadratic extensions. It is believed that noncommutative tori with real multiplication may play an analogous role in the study of real quadratic extensions of $\mathbb{Q}$. In order to achieve arithmetical applications it is important to realize algebraic structures underlying noncommutative tori. This is our main motivation for the study of the homogeneous coordinate rings described below.


\subsection{Complex structures on tori and holomorphic connections}

A complex structure on the noncommutative torus $\mathbb{T}_\theta^2$ is determined through the derivations $\delta_{1}$ and $\delta_{2}$ on $\mathcal{A}_\theta$ by choosing a complex structure on the Lie algebra 
$L = \mathbb{R}^{2}$ of $\mathbb{T}^2$. For this we make a decomposition of the complexification of $L$ into two complex conjugate subspaces. This can be done by choosing a complex parameter $\tau$ with nonzero imaginary part and taking $\{ 1, \tau \}$ as a basis for the holomorphic part of this decomposition. The resulting derivation $\delta_{\tau} = \tau \delta_1 + \delta_{2}$ is by definition a complex structure on $\mathbb{T}_\theta^2$. 
Explicitly we have 
\begin{eqnarray}
\delta_\tau: \sum_{n,m \in \mathbb{Z}} a_{n,m}U^nV^m \mapsto (2 \pi \imath )\sum_{n,m \in \mathbb{Z}} (n\tau + m)a_{n,m}U^nV^m 
\end{eqnarray}
This derivation should be viewed as an analog of the operator $\bar{\partial}$ on a complex elliptic curve. We will denote by $\mathbb{T}_{\theta, \tau}^2$ the noncommutative torus $\mathbb{T}_\theta^2$ equipped with this complex structure. In what follow we will assume that $\mathrm{Im}(\tau) < 0$. We will also freely refer to $\tau$ as the complex structure on $\mathbb{T}_{\theta, \tau}^2$.

Complex structures on noncommutative tori were introduced by Connes in relation with the Yang-Mills equation and positivity in Hochschild cohomology for noncommutative tori (cf. 
\cite{Connes2}). The study of the structure of the space of connections associated to the above derivations was carried out in \cite{ConnesRieffel}. An approach through noncommutative analogs of theta functions was developed in \cite{Schwarz1,DiS}, where these are viewed as holomorphic sections of vector bundles on noncommutative tori. The resulting categories were studied throughly in \cite{PoS} and \cite{Polishchuk2}.
A holomorphic structure on a  right $\mathcal{A}_\theta$-module $E$ is given by an operator $\bar{\nabla}:E\rightarrow E$ which is compatible with the complex structure $\delta_\tau$ in the sense that it satisfies the following Leibniz rule:
\begin{eqnarray}
\bar{\nabla}(ea)= \bar{\nabla}(e)a + e \delta_\tau (a), \quad e \in
E, a\in\mathcal{A}_\theta 
\end{eqnarray}

Given a holomorphic structure $\bar{\nabla}$ on a right $\mathcal{A}_\theta$-module $E$ the corresponding set of holomorphic sections is the space 
\begin{eqnarray}
H^0(\mathbb{T}_{\theta,\tau}^2, E_{\bar{\nabla}}):= \mathrm{Ker} (\bar{\nabla})
\end{eqnarray}

On every basic module $E_{d,c}$ one can define a family of holomorphic structures $\{ \bar{\nabla}_z \}$ depending on a complex parameter $z\in \mathbb{C}$:
\begin{eqnarray}
\label{nabla}
\bar{\nabla}_z(f)=\frac{\partial f}{\partial x} + 2 \pi \imath \left( \frac{ d \tau}{c \theta + d}
x+z \right )f.
\end{eqnarray}

By definition a {\it standard holomorphic vector bundle on $\mathbb{T}_{\theta, \tau}^2$} is given by a basic module $E_{d,c}=E_{g}$ together with one of the holomorphic structures $\bar{\nabla}_z$.

The spaces of holomorphic sections of a standard holomorphic vector bundle on $\mathbb{T}_{\theta, \tau}^2$ are finite dimensional (cf. \cite{PoS}, Section 2). If $c\theta + d > 0$ then $\dim H^0(E_{g}, \bar{\nabla}_0 ) = c$. On what follows we will consider the spaces of holomorphic sections corresponding to  $ \bar{\nabla}_0$:
\begin{eqnarray}
\label{H0}
\mathcal{H}_g := H^0(\mathbb{T}_{\theta,
\tau}^2, E_{g, \bar{\nabla}_0}).
\end{eqnarray}

A basis of $\mathcal{H}_g$ is given by the Schwartz functions:
\begin{eqnarray}
\label{bases1}
\varphi_{\alpha}(x,\beta) = \exp (-\frac{c \tau }{c\theta +d} \frac{x^2}{2})\delta^{\beta}_{\alpha}  \qquad \alpha =1,...,c.
\end{eqnarray}

The tensor product of holomorphic sections is again holomorphic. Using the above basis the product can be written in terms of the corresponding structure constants.

\begin{theorem}\emph{(\cite{PoS} Section 2)}
\label{poli1}
Suppose $g_{1}$ and $g_{2}$ have positive degree. Then $g_1g_2$ has positive degree and 
$t_{g_{1}, g_{2}}$ induces a well-defined linear map
\begin{eqnarray}
\label{producto2}
t_{g_{1}, g_{2}}: \mathcal{H}_{g_{1}}(g_{2}\theta)\otimes_{\mathbb{C}}\mathcal{H}_{g_{2}} (\theta) \rightarrow \mathcal{H}_{g_{1} g_{2}}(\theta).
\end{eqnarray}

Let $g_1,g_2$ and $g_1g_2$ be given by  
$$
g_1= \left(
\begin{array}{cc}
  a_1 & b_1 \\
  c_1 & d_1 \\
\end{array}
\right), \quad 
g_2= \left(
\begin{array}{cc}
  a_2 & b_2 \\
  c_2 & d_2 \\
\end{array}
\right), \quad 
g_1 g_2 = \left(
\begin{array}{cc}
  a_{12} & b_{12} \\
  c_{12} & d_{12} \\
\end{array}
\right)
$$
and let $\{ \varphi_{\alpha} \}$, $\{\varphi_{\beta}' \}$ and $ \{ \psi_{\gamma} \}$ be respectively the bases of $\mathcal{H}_{g_{1}}(g_{2}\theta)$, $\mathcal{H}_{g_{2}} (\theta)$ and $\mathcal{H}_{g_{1} g_{2}}(\theta)$ as given in (\ref{bases1}). Then 
\begin{eqnarray}
\label{stru1}
t_{g_{1}, g_{2}}:  \varphi_{\alpha}\otimes \varphi_{\beta}' \mapsto C^{\gamma}_{\alpha , \beta} \psi_{\gamma}
\end{eqnarray}

\noindent Where for $\alpha= 1,...,c_1$, $\beta = 1,...,c_2$ and $\gamma = 1,..., c_{12}$ we have
\begin{eqnarray}
\label{C1}
C^{\gamma}_{\alpha , \beta} = \sum_{m\in I_{g_1, g_2}(\alpha , \beta, \gamma)} 
\exp [-\pi \imath \frac{ \tau m^{2}}{2  c_1  c_2 c_{12}} ]
\end{eqnarray}
with
\begin{eqnarray*}
I_{g_1, g_2}(\alpha , \beta, \gamma) = \{ n\in \mathbb{Z} &|&  n\equiv -c_1\gamma + c_{12}\alpha \mod  c_{12}c_1, \; \\
&&n\equiv c_{2}d_{12} \gamma - c_{12}d_2 \beta \mod c_{12}c_2 \}
\end{eqnarray*}
\end{theorem}
(we use the convention of summing over  repeated indices).

\subsection{Homogeneous coordinate rings}

Given a projective scheme $Y$ over a field $k$ together with an ample line bundle 
$\mathcal{L}$ on $Y$ one can construct the homogeneous coordinate ring
\begin{eqnarray*}
B = \bigoplus_{n\geq 0} H^0(Y, \mathcal{L}^{\otimes n}).
\end{eqnarray*}
This ring plays a prominent role in the study of the geometry of $Y$.

In \cite{Po} Polishchuk proposed an analogous definition of the homogeneous coordinate ring of a real multiplication noncommutative torus $\mathbb{T}_{\theta}^2$ in terms of holomorphic sections of tensor powers of a standard holomorphic vector bundle on $\mathbb{T}_{\theta, \tau}^2$. As mentioned above the real multiplication condition is fundamental in order to be able to perform the tensor power operation on holomorphic bundles over 
$\mathbb{T}_{\theta, \tau}^2$.

Assume that $\theta \in \mathbb{R}$ is a quadratic irrationality. So there exists some $g \in SL_2(\mathbb{Z})$ with $g\theta= \theta$ and $\mathbb{T}_{\theta}^2$ has real multiplication. Fix a complex structure $\tau$ on $\mathbb{T}_{\theta}^2$. In the case 
$E=E_g(\theta)$ we can extend a holomorphic structure on $E_g$ to a holomorphic structure on the tensor powers $E_{g}^{\otimes n}$. Following \cite{Po} we define a {\it homogeneous coordinate ring for $\mathbb{T}_{\theta,\tau}^2$} by
\begin{eqnarray}
\label{LAdefinicion}
B_g(\theta,\tau) &=& 
\bigoplus_{n\geq 0} H^0(\mathbb{T}_{\theta,\tau}^2, E_{\bar{\nabla}_0}^{\otimes n})\\
 &=& \bigoplus_{n\geq 0} \mathcal{H}_{g^{n}} \nonumber
\end{eqnarray}

The category of holomorphic vector bundles on $\mathbb{T}_{\theta, \tau}^2$ is equivalent to the heart $\mathcal{C}^{\theta}$ of a t-structure on
$\mathcal{D}^{b}(E_{\tau})$, the derived category of the elliptic curve 
$E_{\tau} = \mathbb{C}/ (\mathbb{Z}\oplus\tau \mathbb{Z})$. In \cite{Po} Polishchuk 
exploits this equivalence in order to study the properties of the algebra 
$B_g(\theta,\tau)$ by studying the the corresponding image under this equivalence (cf. \cite{Mahanta}).

The following result characterizes some structural properties of $B_g(\theta,\tau)$ in terms of the matrix elements of $g$:
\begin{theorem}\emph{(\cite{Po} Theorem 3.5)}
\label{poli2}
Assume $g\in SL_2 (\mathbb{Z})$ has positive real eigenvalues: 
\begin{enumerate}
\item If $c \geq a+ d$ then $B_g(\theta,\tau)$ is generated over $\mathbb{C}$ by $\mathcal{H}_{g}$.
\item If $c \geq a+ d+1$ then $B_g(\theta,\tau)$ is a quadratic algebra.
\item \label{la3} If $c \geq a+ d+2$ then $B_g(\theta,\tau)$ is a Koszul algebra.
\end{enumerate}
\end{theorem}

Let us recall these definitions. If $A=\bigoplus_{n\geq 0}A_n$ is a connected graded algebra over a field $k$ generated by its degree one piece 
$A_1$ then $A$ is isomorphic to a quotient 
$T(A_1)/ \mathcal{I}$ where $T(A_1)= \bigoplus_{n\geq 0}A_1^{\otimes n} $ is the tensor algebra of the vector space $A_1$ and $\mathcal{I}$ is a two-sided ideal in 
$T(A_1)$. The algebra $A$ is a \emph{quadratic algebra} if the ideal $\mathcal{I}$ can be generated by homogeneous elements of degree two. Since 
$A$ is connected we can consider $A_0=k$ as a left module over $A$. A quadratic algebra $A$ is a \emph{Koszul algebra} if the graded $k$-algebra $\bigoplus_{n\geq 0} \mathrm{Ext}_{A}^n(k,k)$ is generated by $\mathrm{Ext}_{A}^1(k,k)\simeq A_1^{*}$. 

In a different but related perspective homogeneous coordinate rings on noncommutative tori have been recently studied by Vlasenko in \cite{Masha} where the theory of rings of quantum theta functions is developed.

\subsection{Arithmetic structures}

As in the last section let $\theta$ be a real quadratic irrationality, assume 
$g \in SL_2(\mathbb{Z})$ is a fractional linear transformation fixing $\theta$ and let $\tau$ be a complex structure on $\mathbb{T}_{\theta}^2$. The restriction to $\mathcal{H}_{g} \otimes \mathcal{H}_{g}$ of the multiplication map in the homogeneous coordinate ring $B_g(\theta,\tau)$ is given by the map $t_{g, g}$ defined in (\ref{stru1}). 
The corresponding structure constants $C^{\gamma}_{\alpha , \beta}$ are of the form
\begin{eqnarray}
\vartheta_{r}(l \tau)
\end{eqnarray}
where $\vartheta_{r}(l \tau)$ is the theta constant with rational characteristic $r$
defined by the series 
\begin{eqnarray}
\label{thetaseries}
\vartheta_{r}(l \tau) = \sum_{n \in \mathbb{Z}} \exp [\pi \imath (n+r)^{2} l \tau]
\end{eqnarray}
and $r \in \mathbb{Q}$, $l\in \mathbb{N}$ depend on $\gamma, \alpha , \beta$ and $g$.
(the reader may refer to \cite{tata1} and \cite{tata3} for the main results about theta constants). The arithmetic consequences of this fact were studied in \cite{yo}. In this section we survey some of these results. Throughout the section we assume $g$ satisfies condition 
(\ref{la3})  of Theorem~\ref{poli2} and so  $B_g(\theta,\tau)$ is a quadratic Koszul algebra.

First we relate the field of definition the quadratic algebra $B_g(\theta,\tau)$ with the field of definition of the elliptic curve with period lattice generated by the complex structure 
$\tau$ on $\mathbb{T}_{\theta}^2$:

\begin{theorem} \emph{(\cite{yo})}
\label{superteorema}
Let $E_{\tau}$ be the elliptic curve $\mathbb{C}/ (\mathbb{Z} \oplus\tau \mathbb{Z})$. Let $k$ be its minimal field of definition. Then the algebra $B_{g}(\theta,\tau)$ admits a rational presentation over a finite algebraic extension $K$ of $k$.
\end{theorem}

The minimal field of definition $k$ of the elliptic curve $E_{\tau}$ is given by adjoining to 
$\mathbb{Q}$ the value of the absolute invariant $j(E_{\tau})$ of the curve $E_{\tau}$.
By (\ref{bases1}) there is an isomorphism $\mathcal{H}_g \simeq \sum_{i=1}^{c}\mathbb{C}x_{i}$; then by Theorem~\ref{superteorema} there is a finite algebraic extension $K$ of $k$ and an embedding $K \hookrightarrow \mathbb{C}$ such that the algebra $B_g(\theta,\tau)$ admits a presentation 
\begin{eqnarray}
B_g(\tau,\theta) =  \mathbb{C}\langle x_1,...,x_c \rangle / \mathcal{R}
\end{eqnarray}
with $\mathcal{R}$ a two sided ideal in the free algebra $\mathbb{C}\langle x_1,...,x_c \rangle$
generated by elements in $(\sum_{i=1}^{c}K x_{i})^{\otimes 2}$. 
In particular it makes sense to restrict scalars and consider the $K$-algebra.
\begin{eqnarray}
B_g(\tau,\theta)_{K} =  K\langle x_1,...,x_c \rangle / \mathcal{R}
\end{eqnarray}

Starting with  Theorem~\ref{superteorema} we can study the properties 
of $B_g(\tau,\theta)$ for special values of $\tau$ giving interesting fields of definition. For instance, one can start with an elliptic curve defined over a number field; in this case $j(E_{\tau})$ is algebraic over $\mathbb{Q}$ and the field $K$ is again a number field.

The fact that the homogeneous coordinate ring $B_g(\tau,\theta)$ can be described in terms of theta constants also has some nice consequences when considering $B_g(\tau,\theta)$ as a family of algebras parametrized by $\tau$. Theta constant with rational characteristics are half-weight modular forms. From this fact a presentation of the algebra 
$B_g(\tau,\theta)$ can be given in such a way that each one of the coefficients in the defining relations is a modular function of certain weight and level. Being more explicit:

\begin{theorem} \emph{(\cite{yo})}
A presentation 
\begin{eqnarray}
B_g(\tau,\theta) =  \mathbb{C}\langle x_1,...,x_c \rangle / \mathcal{R}
\end{eqnarray} 
can be given in such a way that the two sided-ideal $\mathcal{R}$ is generated by elements of the form: 
\begin{eqnarray*}
\mathit{f} = {v}^{1} x_{i_1} x_{j_1} +...+ {v}^{c} x_{i_c} x_{j_c}
\end{eqnarray*}
where, as a function of $\tau$, each ${v}^{k}$ is a modular form of cusp type.
\end{theorem}
This result can be used in order to define quadratic algebras associated to 
$\mathbb{T}_{\theta}^2$ independently of the choice of a particular complex structure $\tau$. Modular forms of cusp type define cohomology classes in universal covers of modular curves and therefore can be paired with canonical cycles on these covers. By averaging the coefficients of the defining relations over these cycles one obtains algebras which correspond to averages of the homogeneous coordinate rings $B_g(\tau,\theta)$ along curves in the space 
$\mathbb{H} = \{ \tau \in \mathbb{C} | \mathrm{Im}(\tau ) < 0 \} $ parameterizing complex structures on 
$\mathbb{T}_{\theta}^2$. Of particular interest are those algebras coming from averages along geodesics in $\mathbb{H}$ connecting points in the set of cusps $\mathbb{Q} \cup \{ \infty \}$. These geodesics correspond to homology classes known as modular symbols 
(cf. \cite{Manin4, Shokurov1}). Considering geodesics in $\mathbb{H}$ with limits in 
$\mathbb{R} \setminus \mathbb{Q}$ the theory of modular symbols can be extended in order to define homology classes corresponding to asymptotic limits of modular symbols (cf. \cite{ManinMarcolli}). In particular Yu. Manin and M. Marcolli have shown that the limiting modular symbol associated to the quadratic irrationality $\theta$ exists and can be computed as a linear combinations of classical modular symbols corresponding to cusps in 
$\mathbb{Q} \cup \{ \infty \}$. Applying the averaging process described above along this linear combination of modular symbols one may obtain a quadratic algebra $B_g(\theta)$ canonically associated to $g$ and $\theta$.


\section{Spherical manifolds}

The fact that the geometry of a Riemannian manifold is encoded in terms of spectral data plays a central role in noncommutative geometry (cf. \cite{Connes4, Connes3, ConnesMarcolli2}). 
This spectral point of view was exemplified by A. Connes in the survey article \cite{Connes4} through the study and characterization of spheres. The non-triviality of the pairing between the $K$-theory of the spheres and the $K$-homology class of the Dirac operator is given in terms of a simple polynomial equation. This non-triviality is fundamental in spectral geometry. The formula of this pairing owes its simple form to the fact that lower terms of the Chern character of canonical elements in the $K$-theory of spheres vanish. The vanishing of these lower terms encode to a large extent the algebraic properties of the algebras of functions over the sphere. 
One remarkable fact is that there are noncommutative spaces that share the geometric and algebraic characteristics exhibited by spheres. 

\subsection{Chern characters}

We start by recalling the definition of the Chern character in $K$-theory. Let $\mathcal{A}$ be a unital $*$-algebra; we denote by $\tilde{\mathcal{A}}$ the quotient of 
$\mathcal{A}$ by the subspace $\mathbb{C} \mathbf{1}_{A}$. Recall that $K_{0}(\mathcal{A})$, 
the even $K$-theory of $\mathcal{A}$, can be defined in terms of equivalence classes of self-adjoint idempotents in the matrix algebras $\textsl{M}_q(\mathcal{A})$. Recall also that $K_{1}(\mathcal{A})$, the odd $K$-theory group of $\mathcal{A}$, can be defined in terms of equivalence classes of unitaries in the matrix algebras $\textsl{M}_q(\mathcal{A})$ (cf. \cite{Masoud}). As before we use the convention of summing over repeated indices.

\begin{definition}
The (even) Chern character of a self-adjoint idempotent 
$e=e^{2}=e^{*} \in \textsl{M}_q(\mathcal{A})$, 
\begin{eqnarray}
\mathrm{ch}(e) &=& \mathrm{ch}_0 (e)+ \mathrm{ch}_1(e) + \mathrm{ch}_2(e) + \dots 
\end{eqnarray}
is given by the homogeneous components $\mathrm{ch}_k(e)\in \mathcal{A} \otimes (\tilde{\mathcal{A}}^{\otimes 2k})$ defined by
\begin{eqnarray}
\mathrm{ch}_k(e) &=& \mathrm{Tr} \left( (e- \frac{1}{2}) \underbrace{ e \otimes \dots \otimes e }_{2k}  \right)   \\ 
\nonumber
&=&  \left(e^{i_0}_{i_1} -\frac{1}{2}\delta^{i_0}_{i_1}
\right)\otimes e^{i_1}_{i_2}\otimes e^{i_2}_{i_3}\otimes \dots \otimes
e^{i_k}_{i_0}
\end{eqnarray}
\end{definition}

\begin{definition}
The (odd) Chern character of a unitary matrix $U\in \textsl{M}_q(\mathcal{A})$, $U^*U=UU^*=\mathbf{1}_{A}$,
\begin{eqnarray}
\mathrm{ch}(U) &=& \mathrm{ch}_{\frac{1}{2}}(U)+ \mathrm{ch}_{\frac{3}{2}}(U) + \mathrm{ch}_{\frac{5}{2}}(U) + \dots 
\end{eqnarray}
is given by the homogeneous components
$\mathrm{ch}_{k+\frac{1}{2}}(U)\in \mathcal{A} \otimes (\tilde{\mathcal{A}}^{\otimes (2k+1)})$
defined by
\begin{eqnarray}
\mathrm{ch}_{k+\frac{1}{2}}(U) &=&  \mathrm{Tr} ( \underbrace{(U-1) \otimes (U^{*}-1)  \otimes \dots \otimes  (U-1) \otimes (U^{*}-1) }_{2k +2}  )   \\ 
\nonumber
&=& U^{i_0}_{i_1}\otimes U^{* i_1}_{i_2}\otimes \dots \otimes U^{* i_{2k+1}}_{i_0} -
U^{*i_0}_{i_1}\otimes U^{ i_1}_{i_2}\otimes
\dots \otimes U^{ i_{2k+1}}_{i_0}
\end{eqnarray}
\end{definition}
(we use the convention of summing over repeated indices).

The above Chern characters define maps between the $K$-theory groups of the algebra 
$\mathcal{A}$ and its cyclic homology groups (cf. \cite{Connes2, Connes4}). In the commutative case (corresponding to the algebra of functions on a smooth manifold) one recovers the classical Chern character from the Atiyah-Hirzebruch $K$-theory to the de Rham homology of the manifold. It is also possible to define a Chern character in $K$-homology. The pairing between the $K$-theory of the algebra $\mathcal{A}$ and its $K$-homology can be computed via the corresponding Chern characters in terms of the pairing between the cyclic homology of the algebra and its cyclic cohomology. Given a spectral triple $(\mathcal{A}, \mathcal{H}, D)$ the operator $D$ defines an element in the $K$-homology of the algebra $\mathcal{A}$. This class plays the role of the fundamental class of the geometry encoded by the triple $(\mathcal{A}, \mathcal{H}, D)$ (cf. 
\cite{ConnesMarcolli2}).

\subsection{The two sphere $S^{2}$} (\cite{Connes4})
Consider the following projector over the sphere $S^2$:
\begin{eqnarray}
e \, = \, \frac{1}{2} \left(
                              \begin{array}{cc}
                                1+z & x- \imath y  \\
                                x+ \imath y & 1 - z \\
                              \end{array}
                            \right) \in \textsl{M}_2(
                            C^{\infty}(S^2))
\end{eqnarray}
where $x,y$ and $z$ are the coordinate functions of $S^2 \subset \mathbb{R}^3$  (so $x^2+y^2+z^2=1$). One checks that $e=e^2=e^*$ and so $e$ is actually a projection. Also
$\mathrm{ch}_{0}(e)=0$ and the element $\mathrm{ch}_{1}(e)= \frac{\imath }{2} [x\otimes (y\otimes z- z\otimes y) - y\otimes (x\otimes z - z\otimes x) + z \otimes (x\otimes y - x\otimes y)] $ corresponds (up to the constant $\frac{\imath}{2}$) to the volume form of the round metric on $S^2$.

One can go the other way around and note that the algebra $C^{\infty}(S^2)$ is generated by the matrix elements of the projection $e$. Moreover, it can be shown that the conditions $e=e^2=e^*$, $\mathrm{ch}_{0}(e)=0$ and $\mathrm{ch}_{1}(e)\neq 0$ determine completely the algebra $C^{\infty}(S^2)$ (cf. 
\cite{Varilly}). 
Once we impose this conditions on a on a $2\times 2$ matrix with entries in a $*$-algebra it follows that such matrix can be written in the form 
\begin{eqnarray}
e \quad = \quad \frac{1}{2} \left(
                              \begin{array}{cc}
                                1+z & x- \imath y  \\
                                x+\imath y & 1-z \\
                              \end{array}
                            \right) 
\end{eqnarray}
where $x,y$ and $z$ are mutually commuting self-adjoint elements of the algebra and satisfy the relation $x^2+y^2+z^2=1$. The joint spectrum of the $*$-algebra generated by $x,y$ and $z$  in any faithful representation is $S^2$.

\subsection{The noncommutative spheres $S^{4}_{\theta}$}(\cite{LandiConnes})

In dimension four interesting phenomena appear. This was first noticed by studying projections satisfying the same $K$-theoretic equations as the corresponding analogs of the above projection for $S^4$. These equations impose relations on the algebras generated by their matrix elements. A. Connes and G. Landi studied the geometry of the corresponding noncommutative spaces as part of their work on noncommutative instantons. Below we recall their construction of noncommutative analogs of $S^4$.

Let $\theta$ be real number and let $C_{alg}(S^{4}_\theta)$ be the unital $*$-algebra generated by 3 elements $a,b$ and $x$ together with the relations
\begin{eqnarray*}
x&=&x^*\\
a^*a &=& aa^* \\
b^*b &=& bb^* \\
ab &=& e^{2\pi \imath \theta}ba \\
a^*b&=&e^{-2\pi \imath \theta}ba^* \\
a^*a&+&bb^*+x^2 \, = \, 1.
\end{eqnarray*}
We call $C_{alg}(S^{4}_\theta)$ the algebra of polynomial functions on the noncommutative 4-sphere $S^{4}_\theta $.

\begin{theorem}\emph{(\cite{LandiConnes})}
Consider the matrix
\begin{eqnarray}
e = \left(
      \begin{array}{cccc}
        1+x & 0  &     a        &   b \\
        0 & 1+x  & -\lambda b^* &  a^* \\
        a^* & -\bar{\lambda} b & 1-x  & 0 \\
        b^* & a & 0 & 1-x  \\
      \end{array}
    \right)  \in \textsl{M}_4(C_{alg}(S^{4}_\theta)), \qquad \lambda =e^{2\pi \imath \theta}
\end{eqnarray}
Then
\begin{itemize}
  \item $e=e^2=e^*$.
  \item $\mathrm{ch}_{0}(e)=\mathrm{ch}_{1}(e)=0$
  \item $\mathrm{ch}_{2}(e)\neq 0$
\end{itemize}
\end{theorem}

The algebra $C_{alg}(S^{4}_\theta)$ can be completed to a pre-$C^{*}$-algebra $C^\infty (S^{4}_\theta)$, this algebra acts on the Hilbert space $\mathcal{H}$ of $L^2$ spinors over on $S^{4}$ and together with the corresponding Dirac operator $D$ forms an even spectral triple $(C^\infty (S^{4}_\theta), \mathcal{H}, D)$ that satisfies all the axioms of a noncommutative spin geometry. 

More generally one can consider a compact manifold $M$ whose automorphism group has rank 2. In this case the torus $\mathbb{T}^2$ acts on the algebra $C^\infty (M)$ and we may twist the product by considering also the action of $\mathbb{T}^2 $ on the noncommutative torus $\mathbb{T}^2_\theta$. More precisely, one defines the pre-$C^{*}$-algebra $C^\infty (M_\theta)$ as the sub-algebra of the product $C^\infty (M) \hat{\otimes} C^\infty (\mathbb{T}^2_\theta)$ consisting of elements which are invariant under the diagonal action of the group $\mathbb{T}^2$. We consider then the  algebra $C^\infty (M_\theta)$ as the algebra of smooth function on the noncommutative manifold 
$M_\theta$. We will refer to the class of noncommutative manifolds obtained in this way and their higher dimensional analogs as \emph{$\theta$-deformations} of classical manifolds (cf. \cite{LandiConnes, CoDVI}).

If the manifold $M$ is a spin manifold and $(C^\infty (M), \mathcal{H}, D)$ is the corresponding commutative spectral triple then $C^\infty (M_\theta)$ acts on the Hilbert space $\mathcal{H}$ and it can be shown  that $(C^\infty (M_\theta),
\mathcal{H}, D)$ is a spectral triple. These kinds of geometries are called \emph{isospectral deformations}.

\subsection{Noncommutative 3-Spheres}

In this section we survey the results of \cite{CoDVII,CoDVIII} where a complete classification of the algebras corresponding to the analog in dimension $3$ of the above situation was given. The study of the moduli space that parameterizes these noncommutative three-dimensional spheres was also carried out in detail in these articles. The tools developed by A. Connes and M. Dubois-Violette for the study of the geometric properties of these noncommutative spaces provide a bridge between the purely algebraic framework of quadratic algebras and the analytic framework of 
$C^{*}$-algebras. The notion of central quadratic form plays a fundamental role in order to control the behavior of the corresponding norms.

By definition \emph{the algebra of polynomial functions on a noncommutative spherical manifold of dimension $3$} is a complex unital $*$-algebra $A$ generated by the elements of a unitary matrix $U\in \mathcal{M}_2(A)$ satisfying 
$\mathrm{ch}_{\frac{1}{2}}(U)=0$ and $\mathrm{ch}_{\frac{3}{2}}(U)\neq 0$. One first looks at the restrictions imposed on the structure of the algebra by the fact that $U$ is a unitary and 
$\mathrm{ch}_{\frac{1}{2}}(U)=0$. It will turn out that the non-vanishing of 
$\mathrm{ch}_{\frac{1}{2}}(U)$ follows from these conditions. Thus, being more explicit, $A$ is a 
unital $*$-algebra generated by the elements of a matrix $U\in \mathcal{M}_2(A)$ satisfying 
\begin{eqnarray}
\label{esp1}
UU^* =  U^*U =1 \\
\label{esp2}
U^{i}_{j} \otimes U^{*j}_{i} -   U^{*i}_{j} \otimes U^{j}_{i} = 0
\end{eqnarray}
As before we use the convention of summing over repeated indices.

It is convenient to consider also the corresponding homogeneous problem coming from relaxing a little the unitarity condition and considering ``unitaries up to scale". This corresponds to studying the algebra $B$ of polynomial functions on a noncommutative plane of dimension $4$ spanned by a spherical manifold of 
dimension $3$. Explicitly $B$ is a unital $*$-algebra generated by the elements of a matrix $U\in \mathcal{M}_2(B)$ satisfying 
\begin{eqnarray}
\label{pl1} 
UU^* = U^*U  \, \in \mathbf{1}_2 \otimes B \\
\label{pl2}
U^{i}_{j} \otimes U^{*j}_{i} - U^{*i}_{j} \otimes U^{j}_{i} =  0
\end{eqnarray} 
where $\mathbf{1}_2$ is the identity matrix in $\mathcal{M}_2(\mathbb{C})$.

Once we have chosen a basis $\{ \tau_{\mu} | \; \mu = 0,\dots , 3 \}$ for 
$\mathcal{M}_2(\mathbb{C})$ we can write any matrix $U\in \mathcal{M}_2(A)$ (resp. $\mathcal{M}_2(B)$) as
\begin{eqnarray} 
U = \tau_{\mu} z^{\mu}, \quad z^{\mu} \in A \, (resp. \,B)
\end{eqnarray} 
and the above relations become relations between the elements $z^{\mu} \in B$. 
We may for instance take the basis
\begin{eqnarray}
\tau_{0}= \mathbf{1}_2, \quad  \tau_{1}=\imath \sigma_{1}, \quad
\tau_{2}=\imath \sigma_{2},\quad \tau_{3}=\imath \sigma_{3}. 
\end{eqnarray} 
where the $\sigma_{k}$ are the Pauli matrices
\begin{eqnarray}
\sigma_1= \left(
\begin{array}{cc}
  0 & 1 \\
  1 & 0\\
\end{array}
\right), \quad 
\sigma_2= \left(
\begin{array}{cc}
  0 & -\imath \\
  \imath & 0 \\
\end{array}
\right), \quad 
\sigma_{3} = \left(
\begin{array}{cc}
  1 & 0 \\
  0 & -1\\
\end{array}
\right) 
\end{eqnarray} 
The basis $\{ \tau_{0}, \tau_{1}, \tau_{2}, \tau_{3} \}$ is orthonormal for the inner product $\langle a , b \rangle =\frac{1}{2}\mathrm{Trace}(a^{*}b)$ in 
$\mathcal{M}_2(\mathbb{C})$.
If we now write $U = \tau_{\mu} z^{\mu}$ the equality (\ref{pl2}) is equivalent to:
\begin{eqnarray}
z^{\mu *} &=& \Lambda_{\nu}^{\mu}z^{\nu} 
\end{eqnarray} 
where $\Lambda \in \mathcal{M}_2(\mathbb{C})$ is a unitary symmetric matrix.
By imposing the conditions of unitarity and unitarity up to scale we are led to the following definitions (cf. \cite{CoDVII,CoDVIII}):
\begin{definition}
Let $\Lambda \in \mathcal{M}_4(\mathbb{C})$ be a unitary symmetric matrix. We
define $C_{alg}(\mathbb{R}^{4}_{\Lambda})$, \emph{the algebra of polynomial functions on the noncommutative four-plane $\mathbb{R}^{4}_{\Lambda}$} as the 
$*$-algebra generated by four elements $\{ z^{0}, z^{1}, z^{2}, z^{3} \}$ subject to the relations
\begin{eqnarray}
\label{rela1}
z^{\mu *} &=& \Lambda_{\nu}^{\mu}z^{\nu} \\
\label{rela2}
z^{k} z^{0 *} -  z^{ 0}z^{k *} &-& \sum \epsilon_{klm} z^{l} z^{m *} \;=\; 0 \\
\label{rela3}
z^{0 *} z^{k }- z^{ k *} z^{0 } &-& \sum \epsilon_{klm} z^{l *} z^{m} \;=\; 0 
\end{eqnarray}
where $(k,l, m)$ runs over the cyclic permutations of 
$(1,2,3)$ and $\epsilon_{klm} $ is the totally antisymmetric tensor.

We define $C_{alg}(S^{3}_{\Lambda})$, \emph{the algebra of polynomial functions on the noncommutative three-sphere $S^{3}_{\Lambda}$} as the quotient of $C_{alg}(\mathbb{R}^{4}_{\Lambda})$ by the 
two-sided ideal generated by
\begin{eqnarray}
\label{rela4}
\sum z^{\mu} z^{\mu *} - 1
\end{eqnarray}
\end{definition}

The following theorem shows that $S^{3}_{\Lambda}$ is a spherical manifold of dimension $3$ and any spherical manifold can be embedded in some  $S^{3}_{\Lambda}$.
\begin{theorem}\emph{(\cite{CoDVII,CoDVIII})}
For any unitary symmetric matrix $\Lambda \in \mathcal{M}_4(\mathbb{C})$ the unitary 
\begin{eqnarray}
U = \tau_{\mu} z^{\mu} \quad  \in  \mathcal{M}_2(C_{alg}(S^{3}_{\Lambda}))
\end{eqnarray} 
satisfies  $\mathrm{ch}_{\frac{1}{2}}(U)=0$

Conversely, for any complex unital $*$-algebra $A$ generated by the elements of a unitary matrix $\tilde{U} \in \mathcal{M}_2(A)$ satisfying $\mathrm{ch}_{\frac{1}{2}}(\tilde{U})=0$  
there exist a unitary symmetric matrix $\Lambda \in \mathcal{M}_4(\mathbb{C})$ and a 
surjective morphism 
\begin{eqnarray}
C_{alg}(S^{3}_{\Lambda}) \rightarrow A
\end{eqnarray}
carrying  $U$ to $\tilde{U}$.
\end{theorem}

Given $\varphi = (\varphi_1,\varphi_2,\varphi_3 )\in \mathbb{T}^3$ we will denote respectively by 
$\mathbb{R}^{4}_{\varphi}$ and $S^{3}_{\varphi}$ the noncommutative $4$-plane and $3$-sphere corresponding to the matrix
\begin{eqnarray}
\Lambda(\varphi) = \left(
\begin{array}{cccc}
1 & 0 & 0 & 0\\
0 & e^{-2 \pi \imath \varphi_1}& 0 & 0\\
0 & 0 & e^{-2 \pi \imath \varphi_2} & 0\\
0 & 0 & 0 & e^{-2 \pi \imath \varphi_3}
\end{array}
\right)
\end{eqnarray}
It can be shown that for any unitary symmetric matrix $\Lambda \in \mathcal{M}_4(\mathbb{C})$ there exists a point in the torus $\varphi \in \mathbb{T}^3$ such that one has isomorphisms 
$\mathbb{R}^{4}_{\Lambda} \simeq \mathbb{R}^{4}_{\varphi}$ and $S^{3}_{\Lambda} \simeq S^{3}_{\varphi}$ given by isomorphisms of the corresponding algebras. Therefore spherical manifolds in dimension three are parameterized by the points in a $3$-torus $\mathbb{T}^3 = S^1\times S^1\times S^1$.

In order to get Hermitian generators for the 
$*$-algebra $C_{alg}(\mathbb{R}^{4}_{\varphi})$ we take $x^{0}=z^{0}$ and 
$x^{k} = e^{- \pi \imath \varphi_k} z^{k}$, $k=1,2,3$. Thus
\begin{eqnarray}
x^{\mu }=x^{\mu *}, \qquad \mu = 0,\dots,3
\end{eqnarray}
Therefore $C_{alg}(\mathbb{R}^{4}_{\varphi})$ can also be realized as the unital $*$-algebra generated by four Hermitian elements $x^{0},x^{1},x^{2}$ and $x^{3}$ together with the relations
\begin{eqnarray}
\label{ReL1}
\cos(\varphi_k)(x^{0}x^{k}-x^{k}x^{0}) &=& 
\imath \sin(\varphi_l-\varphi_m)(x^{l}x^{m}+x^{m}x^{l}), \\
\label{ReL2}
\cos(\varphi_l-\varphi_m)(x^{l}x^{m}-x^{m}x^{l})&=&
- \imath \sin(\varphi_k)(x^{0}x^{k}+x^{k}x^{0})
\end{eqnarray}
where $(k,l,m)$ runs over the cyclic permutations of $(1,2,3)$.

Note that the value $\varphi = (0,0,0)$ corresponds to the polynomial algebras of the usual 4-plane and 3-sphere.

The first step in the study of these algebras is the identification of a vector field on $\mathbb{T}^3$ corresponding to the foliation coming from the equivalence relation
on $\mathbb{T}^3$ given by
\begin{eqnarray}
\varphi \sim \varphi' & \text{if and only if } \, C_{alg}(\mathbb{R}^{4}_\varphi)\simeq
C_{alg}(\mathbb{R}^{4}_{\varphi'})
\end{eqnarray}
The algebras corresponding to critical values of the foliation will give rise either to  spherical manifolds that are noncommutative $3$-spheres whose suspension is of the form  $S^{4}_\theta$ or to noncommutative $3$-spheres arising from quantum groups. The study of the generic case of the algebras $C_{alg}(S^{3}_{\varphi})$ corresponding to regular values of the foliation uses tools coming from noncommutative algebraic geometry and  nontrivial refinements of these tools in the $C^{*}$-algebra context.

By construction the algebras $C_{alg}(\mathbb{R}^{4}_\varphi)$ are quadratic algebras with 4 generators of degree 1. More explicitly, the algebra $C_{alg}(\mathbb{R}^{4}_\varphi)$ is isomorphic to the quotient 
\begin{eqnarray}
\mathbb{C} \langle x^{0},x^{1},x^{2}, x^{3}  \rangle / \mathcal{R}
\end{eqnarray}
where $\mathbb{C} \langle x^{0},x^{1},x^{2}, x^{3}  \rangle$ is the free $*$-algebra in four Hermitian generators and $\mathcal{R}$ is the two-sided ideal generated by the relations (\ref{ReL1}) and (\ref{ReL2}). 

If we let $V := \sum \mathbb{C} x_i \simeq \mathbb{C}^{4}$ then the algebra 
$ \mathbb{C} \langle x^{0},x^{1},x^{2}, x^{3}  \rangle $ is isomorphic to the tensor algebra of the space $V$:
\begin{eqnarray}
\mathcal{T}(V) = \bigoplus_{n \geq 0}V^{\otimes n} 
\end{eqnarray}
and we can identify the relations (\ref{ReL1}) and (\ref{ReL2}) with elements in 
$V\otimes V$.

To each element $f \in V\otimes V$ we associate a bilinear form 
$\check{f}: V^{*}\times V^{*}\rightarrow \mathbb{C}$. We call $\check{f}$ the 
\emph{multilinearization} of $f$. Since $\check{f}$ is bihomogeneous its zero locus defines a hypersurface in $\mathbb{P}(V)\times\mathbb{P}(V) \simeq \mathbb{P}^{3}(\mathbb{C})\times \mathbb{P}^{3}(\mathbb{C})$. 
The locus of common zeroes of the multilinearizations of the elements of $\mathcal{R}$ defines a variety $\{ \check{f}_i = 0\} = \Gamma\subset \mathbb{P}^{3}\times \mathbb{P}^{3}$. Let $Y_1$ and $Y_2$ be the corresponding projections and $\sigma: Y_1 \rightarrow Y_2$  be the correspondence with graph $\Gamma$. Assume we can make an identification  
$Y_1=Y_2=Y$. If $\sigma$ is an isomorphism we consider it as an automorphism of $Y$. Let $i:Y\hookrightarrow \mathbb{P}^{3}$ be the inclusion and take $\mathcal{L}=i^{*}\mathcal{O}_{\mathbb{P}^{3}}(1)$ where $\mathcal{O}_{\mathbb{P}^{3}}(1)$ is the canonical bundle on $\mathbb{P}^{3}$. The \emph{geometric data associated to the quadratic algebra $C_{alg}(\mathbb{R}^{4}_\varphi)$} is by definition the triple 
\begin{eqnarray}
(Y, \sigma,\mathcal{L} ) .
\end{eqnarray}
The variety $Y$ is called \emph{the characteristic variety $C_{alg}(\mathbb{R}^{4}_\varphi)$}.

Starting from such a triple one can construct the graded algebra:
\begin{eqnarray}
B(Y, \sigma,\mathcal{L} ) = \bigoplus_{n\geq 0} 
H^0(Y,\mathcal{L} \otimes\mathcal{L}^{\sigma}\otimes \dots \otimes\mathcal{L}^{\sigma^{n-1}}) 
\end{eqnarray}
where $\mathcal{L}^{\sigma} := \sigma^{*} \mathcal{L}$ 
and the multiplication of two sections 
$s_{1}\in B(Y, \sigma,\mathcal{L} )_n,s_{2}\in B(Y, \sigma,\mathcal{L} )_m$ is given by $s_1s_2:=s_1\otimes s_2^{\sigma^{n}}$.

The study of the geometric data of $C_{alg}(\mathbb{R}^{4}_\varphi)$ for different values of $\varphi$ is a fundamental step in the classification of these algebras. 
For generic values of the parameter $\varphi$ the characteristic variety $Y$ is given by an elliptic curve $E_{\varphi}$ together with four isolated points and the morphism $\sigma$ acts by translation on the elliptic curve $E_{\varphi}$ and trivially on the four points (cf. \cite{Smith,CoDVII,CoDVIII}).

By construction there is a morphism 
\begin{eqnarray}
C_{alg}(\mathbb{R}^{4}_\varphi) \rightarrow B(Y, \sigma,\mathcal{L} )
\end{eqnarray}

In order to refine this morphism  A. Connes and M. Dubois-Violette introduced in \cite{CoDVII} the notion of a \textit{central quadratic form} on a connected component $C$ of $Y \times Y$. Roughly speaking these are symmetric tensors on the space of generators of the algebra whose action on the points of 
$C$ gives rise to an Hermitian structure on the line bundle $\mathcal{L}$ compatible with the action of $\sigma$ on its sections.

By construction the algebras $C_{alg}(\mathbb{R}^{4}_\mathbf{\varphi})$ are quadratic $*$-algebras, thus; there is an involution preserving the space $V$ of generators
and giving this space a real structure. For this kind of algebras it makes sense to ask whether a compact connected component $C$ of the graph giving the geometric data is invariant under the involution coming from the real structure on $V$ and, this being the case, whether on $C$ the involution commutes with the action of the
automorphism $\sigma$. When these compatibility conditions hold one may obtain a  $*$-homomorphism to a twisted crossed product $C^{*}$-algebra $C(F)\times_{\sigma , \mathcal{L}} \mathbb{Z} $ where $F$ is the image of $C$ under the first projection. In the case of the noncommutative spheres $S^{3}_\mathbf{\varphi}$ corresponding to generic values of $\varphi$ this construction leads to a unital
$*$-morphism
\begin{eqnarray}
C_{alg}(S^{3}_\varphi) \rightarrow C^{\infty}(E_\varphi)\times_{\sigma , \mathcal{L}} \mathbb{Z}
\end{eqnarray}

The existence of this morphisms  is the basis for the study of the geometry of the spheres $S^{3}_\varphi$. In particular using this
map we can pair canonical Hochschild classes in $C^{\infty}(E_\varphi)\times_{\sigma , \mathcal{L}} \mathbb{Z} $ with the corresponding pullback of  $\mathrm{ch}_{\frac{3}{2}}(U)$; since the result obtained is non zero it follows that
\begin{eqnarray}
\mathrm{ch}_{\frac{3}{2}}(U) \neq 0 
\end{eqnarray}
The computation of this pairing is expressed naturally in terms of modular functions, thus providing a bridge between these noncommutative spaces and arithmetic geometry (cf. \cite{CoDVIII} Theorem 12.1). 

\section{Epilogue}

There is a wide variety of quadratic algebras arising naturally from geometric considerations in the framework of algebraic noncommutative geometry 
(cf. \cite{ArtinTate,Smith,StaffordVdBergh, Po}). Applying the tools developed by 
A. Connes and M. Dubois-Violette in \cite{CoDVIII} to the study of these algebras may provide valuable information about the structures which play a role in the interplay between algebraic noncommutative geometry and differential noncommutative geometry.
In particular, applying these tools to the homogeneous coordinate rings $B_g(\tau,\theta)$ and $B_g(\theta)$ associated to real multiplication noncommutative tori could have potential applications. Nontrivial extensions of the theory may have to be developed in order to handle these quadratic algebras since typically they exhibit ill behaved growth properties (cf. \cite{Po}).

Being able to construct embeddings from algebraically defined rings into $C^{*}$-algebras gives important information about their representation theory. In the case of rings coming from arithmetic data  these embeddings may provide crucial information towards applications. This can be seen by considering the recent developments relating class field theory  to quantum statistical mechanics (cf.\cite{BostConnes,ConnesMarcolli,CMR1,CMR2,Benoit,ConsaniMarcolli}). The algebra of observables of a quantum statistical mechanical system is a $C^{*}$-algebra, and the time evolution of the system singles out particular states which are obtained as extremal points in the space of equilibrium states of the system. In quantum statistical mechanical systems describing the explicit class field theory of a particular global field $K$ the extremal equilibrium states are evaluated at observables corresponding to an arithmetically defined rational subalgebra of observables. The values obtained in this way are generators of the abelian extensions of the field $K$. In our context the tools developed in \cite{CoDVIII} could provide a way to construct in a canonical way a quantum statistical mechanical system with a particular given ring as an arithmetic subalgebra. In particular, in view of Manin's real multiplication program outlined above, it seems that investigating the homogeneous coordinate rings of real multiplication noncommutative tori in this context could throw some light about the explicit class field theory problem for real quadratic fields.

\bibliographystyle{amsalpha}

\end{document}